\input amstex
\magnification=\magstep1 \baselineskip=13pt
\documentstyle{amsppt}
\vsize=8.7truein \CenteredTagsOnSplits \NoRunningHeads
\def\today{\ifcase\month\or
  January\or February\or March\or April\or May\or June\or
  July\or August\or September\or October\or November\or December\fi
  \space\number\day, \number\year}

\def\PP{{\bold P}}
\def\EE{{\bold E\thinspace }}
\def\vl{\operatorname{vol}}
\def\per{\operatorname{per}}
\def\Mat{\operatorname{Mat}}

\def\hd{\operatorname{head}}
\def\tl{\operatorname{tail}}
\topmatter
\title Brunn-Minkowski  Inequalities for Contingency Tables 
and Integer Flows \endtitle
\author Alexander Barvinok \endauthor
\address Department of Mathematics, University of Michigan, Ann Arbor,
MI 48109-1043 \endaddress \email barvinok$\@$umich.edu
\endemail
\date March  2006 \enddate
\thanks This research was partially supported by NSF Grant DMS 0400617.
\endthanks
\abstract Given a non-negative $m \times n$ matrix $W=\left(w_{ij}\right)$ and positive integer 
vectors $R=\left(r_1, \ldots, r_m\right)$ and $C=\left(c_1, \ldots, c_n\right)$, 
we consider the total weight $T(R, C; W)$ of $m \times n$ non-negative integer matrices (contingency tables) $D$ with the row sums $r_i$,
the column sums $c_j$, and the weight of $D=\left(d_{ij}\right)$ equal to 
$\prod_{ij} w_{ij}^{d_{ij}}$. In particular, if $W$ is a 0-1 matrix, 
$T(R, C; W)$ is the number of integer feasible flows in a 
bipartite network. We prove a version of the Brunn-Minkowski 
inequality relating the numbers $T(R, C; W)$ and $T(R_k, C_k; W)$, where 
$(R, C)$ is a convex combination of $(R_k, C_k)$ for $k=1, \ldots, p$.
 \endabstract
\endtopmatter
\document

\head 1. Introduction\endhead

\subhead (1.1) The Brunn-Minkowski inequality \endsubhead
The famous Brunn-Minkowski inequality states that for bounded Borel sets 
$A, B \subset {\Bbb R}^d$ and non-negative numbers $\alpha, \beta$ 
such that $\alpha + \beta =1$ one has
$$\vl(\alpha A + \beta B) \geq \vl^{\alpha}(A) \vl^{\beta}(B),$$
where $\vl$ is the usual volume (Lebesgue measure) in Euclidean space ${\Bbb R}^d$
and 
$$\alpha A+ \beta B=\{ \alpha x+ \beta y: \quad x \in A, y \in B \}.$$
The inequality extends to finite families of sets in an obvious way:
if $A_1, \ldots, A_p \subset {\Bbb R}^d$ are bounded Borel sets and 
$\alpha_1, \ldots, \alpha_p$ are non-negative numbers such that 
$\alpha_1 + \ldots + \alpha_p=1$ then 
$$\vl\left(\alpha_1 A_1 + \ldots + \alpha_p A_p\right) \geq \prod_{k=1}^p \vl^{\alpha_k}\left(A_k\right).
\tag1.1.1$$
The Brunn-Minkowski inequality plays an important role in almost all branches of mathematics, 
see \cite{Ga02} for a survey. Inequality (1.1.1) was extended and generalized in numerous 
direction. In particular, we need its functional version, known as the Pr\'ekopa-Leindler
inequality:
\medskip
let  $\alpha_1, \ldots, \alpha_p$ be non-negative numbers such that 
$\alpha_1 + \ldots + \alpha_p=1$ and let 
$g, h_1, \ldots, h_p: {\Bbb R}^d \longrightarrow {\Bbb R}$ be Borel measurable non-negative 
functions such that 
$$g\left(\alpha_1 x_1 + \ldots + \alpha_p x_p \right) \geq \prod_{k=1}^p  h^{\alpha_k}_k(x_k) 
\quad \text{for all} \quad x_1, \ldots, x_k \in {\Bbb R}^d.$$
Then 
$$\int_{{\Bbb R}^d} g(x) \ dx \geq \prod_{k=1}^p \left(\int_{{\Bbb R}^d} h_k(x) \ dx \right)^{\alpha_k},
\tag1.1.2$$ 
 see for example, Section 6.1 of \cite{Vi03} and Section 2.2 of \cite{Le01}. We note 
 that (1.1.1) is obtained from (1.1.2) by choosing $h_k$ to be the indicator function of
 $A_k$, so that $h_k(x)=1$ if $x \in A_k$ and $h_k(x)=0$ if $x \notin A_k$ and $g$ to be 
 the indicator of $\alpha_1 A_1 + \ldots + \alpha_p A_p$. The inequality 
 (1.1.2) remains valid if $dx$ is replaced by a log-concave measure.
\bigskip
In this paper we obtain versions of inequality (1.1.1), respectively (1.1.2), for the number of 
integer points, respectively for the number of weighted integer points, in some special 
polytopes, known as {\it flow polytopes}. 

\subhead (1.2) Contingency tables \endsubhead
Let $R=\left(r_1, \ldots, r_m\right)$ and $C=\left(c_1, \ldots, c_n \right)$ be positive integer 
vectors such that
$$\sum_{i=1}^m r_i = \sum_{j=1}^n c_j=N.$$
 An $m \times n$ non-negative integer matrix $D=\left(d_{ij}\right)$ with 
the row sums $r_1, \ldots, r_m$ and the column sums $c_1, \ldots, c_n$ is called 
a {\it contingency table} with margins $R$ and $C$. Geometrically, one can think of 
the set of contingency tables with prescribed margins as of the set of integer points 
in the {\it transportation polytope} $P(R,C)$ of $m \times n$ matrices $X=\left(x_{ij}\right)$ satisfying 
the equations 
$$\sum_{j=1}^n x_{ij}=r_i \quad \text{for} \quad i=1, \ldots, m, \quad \sum_{i=1}^m x_{ij}=c_j 
\quad \text{for} \quad j=1, \ldots, n$$
and inequalities
$$x_{ij} \geq 0 \quad \text{for all} \quad i,j.$$

The numbers of contingency tables with prescribed margins are of interest because 
of their applications in statistics, combinatorics, and representation theory, see
\cite{DE85}, \cite{DG95}, and \cite{DG04}. 

We consider the number of {\it weighted} tables, defined as follows.
\definition{(1.3) Definition} Let $W=\left(w_{ij}\right)$ be an $m \times n$ non-negative matrix.
For $R=\left(r_1, \ldots, r_m\right)$ and $C=\left(c_1, \ldots, c_n \right)$, we define 
$$T(R,C; W)=\sum_D \prod_{ij} w_{ij}^{d_{ij}},$$
where the sum is taken over all $m \times n$ contingency tables $D=\left(d_{ij}\right)$ 
with the margins $(R, C)$. We agree that $0^0=1$. 
\enddefinition
Geometrically, $T(R,C; W)$ is the generating function over the set of integer points in a 
transportation polytope.  We get the number of points if we choose $W={\bold 1}$, the 
matrix of all 1s.
  
\subhead (1.4) Integer flows \endsubhead
Let $G=(V,E)$ be a directed graph with the set $V$ of vertices, the set $E$ of edges,
without multiple edges or loops. Suppose that an integer $a(v)$, called the {\it excess} 
$v$, is assigned to every vertex $v \in V$ so that 
$$\sum_{v \in V} a(v)=0.$$
A collection $x(e): e \in E$ of non-negative integers is called an {\it integer feasible flow} in 
$G$ if the balance condition is satisfied at every vertex
$$\sum_{e:\  \hd(e)=v} x(e) - \sum_{e:\  \tl(e)=v} x(e)=a(v) \quad \text{for all} \quad v \in V.$$   
If $G$ does not contain directed cycles $v_1 \rightarrow v_2  \rightarrow 
\ldots \rightarrow  v_k \rightarrow v_1$ then the set of feasible flows is compact, so the number of 
 integer feasible flows is finite.

Some interesting quantities can be defined as the number of integer feasible flows in an 
appropriate network. For example, we get the {\it Kostant partition function} 
(for the $A_{n-1}$ root system) if 
$G=K_n$ is a complete graph with the set of vertices $V=\{1, \ldots, n\}$ and edges
$E=\{ i \rightarrow j: \ i>j \}$, cf. \cite{B+04}. Given an integer vector 
$a=\left(a_1, \ldots, a_n\right)$ such that $a_1+ \ldots + a_n=0$, the number $\phi(a)$ 
of integer feasible flows in $K_n$ with the excess at $i$ equal $a_i$ is the value of the 
Kostant partition function at $a$.   

Given a directed graph $G$ on $|V|=n$ vertices and excesses $a(v)$ at its vertices, one can construct
an $n \times n$ matrix $W=\left(w_{ij}\right)$ with $w_{ij} \in \{0,1\}$, a 
vector $R=\left(r_1, \ldots, r_n\right)$ of row sums and a vector $C=\left(c_1, \ldots, c_n\right)$ 
of column sums so that $T(R,C; W)$ is equal to the number of integer feasible flows 
in $G$. To that end, we identify $V=\{1, \ldots, n\}$. Given
the excess $a_i$ at the vertex $i$ of $G$, we find an a priori upper bound   
$z_i \geq 0$ on the total incoming flow to $i$ and let  $r_i=z_i-a_i$ and $c_i=z_i$.
Finally, we let $w_{ij}=1$ if $i=j$ or $i \rightarrow j$ is an edge of $G$ and let 
$w_{ij}=0$ otherwise. 
 
With a feasible flow $\left\{x_e: \ e \in E \right\}$ in $G$, we associate a contingency table 
$D=\left(d_{ij}\right)$ as follows: we let $d_{ij}=x(e)$ provided $i=\hd(e)$ and $j=\tl(e)$ and
let 
$$d_{ii}=r_i-\sum_{e: \ \tl(e)=i} x(e)=c_i-\sum_{e: \ \hd(e)=i} x(e).$$ 
Further, we let $d_{ij}=0$ if $w_{ij}=0$. One can observe that this correspondence is a 
bijection between the integer feasible flows in $G$ and the contingency tables enumerated 
by $T(R,C; W)$.

For example, for the Kostant partition function, we let 
$w_{ij}=1$ if $i \geq j$ and $w_{ij}=0$ otherwise and define $r_1=0$,
$r_i=a_1 + \ldots  + a_{i-1}$ for $i >1$ and  $c_j=a_1 + \ldots  + a_j$ for $j \geq 1$.
Noticing that $r_1=c_n=0$, we cross out the first row and the $n$th column and obtain 
the following description of the Kostant partition function.
\bigskip
Let us define the $(n-1) \times (n-1)$ matrix $W=\left(w_{ij}\right)$ by 
$$w_{ij}=\cases 1 &\text{if\ } i \geq j-1\\ 0 & \text{otherwise.} \endcases$$
Let 
$$r_k=c_k=\sum_{i=1}^k a_k \quad \text{for} \quad k=1, \ldots, n-1.$$
Then  the Kostant partition function $\phi$ satisfies
$$\phi(a_1, \ldots, a_n)=T(R, C; W).$$
 \bigskip
A version of the integer flow enumeration problem involves positive integer {\it capacities} $c(e)$ of edges and requires 
feasible flows to satisfy $x(e) \leq c(e)$. Given a directed 
graph $G$ with capacities one can construct a directed graph $G'$ without capacities so that 
the integer feasible flows in  $G'$ are in a bijection with the integer feasible flows in $G$.
For that, an extra vertex is introduced for every edge of $G$ with capacity, see 
\cite{B+04}. 
  
 \head 2. Main results \endhead

Our main result is the following inequality relating numbers $T(R, C; W)$ of weighted 
contingency tables for different margins $R$ and $C$.
 
 \proclaim{(2.1) Theorem} 
 For a positive integer vector
$B=\left(b_1, \ldots, b_l\right)$
we define 
$$|B|=\sum_{i=1}^l b_i \quad 
 \text{and} \quad \omega(B)=\prod_{i=1}^l {b_i^{b_i} \over b_i!}.$$ 
 
 Let $W=\left(w_{ij}\right)$ be a non-negative $m \times n$ matrix,
 let $R_1, \ldots, R_p$ be positive integer $m$-vectors and let $C_1, \ldots, C_p$ be positive 
 integer $n$-vectors such that 
 $$|R_1| =\ldots =|R_p|=|C_1|=\ldots =|C_p|=N.$$
 Suppose further that $\alpha_1, \ldots, \alpha_p \geq 0$ are numbers such that 
$\alpha_1 + \ldots + \alpha_p=1$.
 Let us define 
 $$R=\sum_{k=1}^p \alpha_k R_k \quad \text{and} \quad C=\sum_{k=1}^p \alpha_k C_k$$
 and suppose that $R$ and $C$ are positive integer vectors.
 
 Then  
 $${N^N \over N!}{T(R,C;W) \over  \omega(R) \omega(C)}  \geq \prod_{k=1}^p
  \left( {T\left(R_k, C_k; W \right) \over \min \bigl\{ \omega(R_k),\  \omega(C_k)\bigr\} } \right)^{\alpha_k}.$$
 \endproclaim
 
 Geometrically, for the transportation polytopes $P(R, C)$ and $P\left(R_k, C_k\right)$ we have 
 $$P(R, C)=\alpha_1 P\left(R_1, C_1\right) + \ldots  + \alpha_pP\left(R_k, C_k\right),$$
 cf. Section 1.2. On the other hand, the corresponding convex combination of integer points 
 in $P\left(R_k, C_k\right)$ does not have to be an integer point in $P(R, C)$. Hence, the 
 existence of an a priori relation between the numbers of integer points in $P\left(R_k, C_k\right)$ and 
 $P\left(R, C\right)$ is not obvious (for a different approach to discrete Brunn-Minkowski 
 inequalities, see \cite{GG01}).
 
What follows is a chain of weaker inequalities which are easier to parse.

\proclaim{(2.2) Corollary} Under the conditions of Theorem 2.1, 
let $$\split &R=\left(r_1, \ldots, r_m\right),\quad C=\left(c_1, \ldots, c_n\right), \quad
a=\min\{m,n\}, \quad \text{and} \\ & s=N/a \quad \text{where} \quad N=\sum_{i=1}^m r_i =
\sum_{j=1}^n c_j. \endsplit$$
Then we have 
\roster
\item
$$ {N^N \over N!} \min\left\{ \prod_{i=1}^m {r_i! \over r_i^{r_i}}, \quad  \prod_{j=1}^n {c_j! \over c_j^{c_j}} \right\} T(R,C; W)
\geq \prod_{k=1}^p T^{\alpha_k} \left(R_k, C_k; W \right).$$
\item $${N^N \over N!} {\Gamma^a(s+1) \over s^N} T(R,C; W)  \geq 
\prod_{k=1}^p T^{\alpha_k} \left(R_k, C_k; W \right).$$ 
\item There is an absolute constant $\kappa>0$ such that
$$(\kappa s)^{{1 \over 2}(a-1)} T(R, C; W) \geq \prod_{k=1}^p T^{\alpha_k} \left(R_k, C_k; W \right).$$ 
\endroster
\endproclaim

Generally speaking, the correction term $(\kappa s)^{(a-1)/2}$ is small compared to 
the value of $T(R, C; W)$. For example, if $w_{ij} \in \{0,1\}$ for all $i, j$ then 
$T(R, C; W)$ is the number of integer points in the flow polytope $P(R, C; W)$ defined 
in the space of $m \times n$ matrices $\left(x_{ij}\right)$ by the equations
$$\split &\sum_{j=1}^n x_{ij}=r_i \quad \text{for} \quad i=1, \ldots, m \\  
&\sum_{i=1}^m x_{ij}=c_j \quad \text{for} \quad j=1, \ldots, n, \quad \text{and} \\
& x_{ij}=0 \quad \text{whenever} \quad w_{ij}=0 \endsplit$$
and inequalities 
$$x_{ij} \geq 0 \quad \text{provided} \quad w_{ij}=1.$$  If we scale $R \longmapsto tR$,
$C \longmapsto tC$ for a positive integer $t$, the number of integer points in 
$P(tR, tC; W)$ grows as a polynomial of $t$ of degree $d=\dim P(R, C; W)$, see, for example, 
Section 4.6 of \cite{St97}, which can be as high as 
$d=(m-1)(n-1)$ in the transportation polytope (see Section 1.2) with $w_{ij} \equiv 1$.
 On the other hand, the correction term $(\kappa s)^{(a-1)/2}$ is 
a polynomial in $t$ of degree $\left(\min\{m,n\}-1\right)/2$.

As another extreme case, let us consider the situation when  the 
numbers $r_i, c_j$ are 
uniformly bounded, while $m$ and $n$ grow. In this case, $T(R,C; W)$ grows roughly as 
$\left(\kappa_1 N\right)^N$,  as long as the number of zeros  in each row and column of 
the 0-1 matrix $W$ is uniformly bounded, cf. \cite{Be74}. The correction term is about $\kappa_2^N$ for 
some absolute constants $\kappa_1, \kappa_2 >0$.

Let us choose an $m \times n$ matrix $c_{ij}$ and let us define matrix $W(t)=\left(w_{ij}(t)\right)$ by
$w_{ij}(t)=\exp\left\{t c_{ij}\right\}$. One can observe that 
$$\lim_{t \longrightarrow +\infty} t^{-1} \ln T\bigl(R,C; W(t)\bigr)=
\max \left\{ \sum_{ij} c_{ij} x_{ij}: \ \left(x_{ij}\right) \in P(R, C) \cap {\Bbb Z}^{m \times n} \right\}.$$
In words: the limit is equal to the maximum value of the linear function defined by matrix 
$\left(c_{ij} \right)$ on the set of integer points in the transportation polytope $P(R, C)$, see Section 1.2. Thus any estimate of the type 
$$\alpha(R, C) T(R, C; W) \geq \prod_{k=1}^p T^{\alpha_k} \left(R_k, C_k; W \right),$$
where $\alpha(R, C)$ is a factor depending on $R$ and $C$ alone, implies that if 
$x_k \in P(R_k, C_k)$ are integer points then the point 
$\alpha_1 x_1 + \ldots + \alpha_k x_k$ lies inside the convex hull of the set of integer points of 
$P(R, C)$, which also follows from the fact that the vertices of $P(R, C)$ are integer.
 
One can ask, naturally, whether the bound in Theorem 2.1 can be strengthened. In particular, the 
following question is of interest:
\bigskip
$\bullet$ Is it true that under conditions of Theorem 2.1, one has
$$T(R, C; W) \geq \prod_{k=1}^p T^{\alpha_k} \left(R_k, C_k; W \right)? \tag2.3$$
Or, perhaps, does the above inequality hold in some interesting special cases, for example, 
when $W={\bold 1}$, the $m \times n$ matrix of all 1s, so that $T(R, C; W)$ is the number 
of contingency tables with the row sums $R$ and column sums $C$? 
\bigskip
There is some circumstantial evidence that the $T(R, C; {\bold 1})$ might indeed satisfy (2.3).
We note that the value of $T(R, C; {\bold 1})$ does not change if the entries of $R$ and 
and $C$ are arbitrarily permuted. Let $a=\left(\alpha_1, \ldots, \alpha_n\right)$ 
and $b=\left(\beta_1, \ldots, \beta_n \right)$ be integer vectors such that 
$$\alpha_1 \geq \alpha_2 \geq \ldots \geq \alpha_n \quad \text{and} \quad \beta_1 \geq \beta_2 \geq \ldots \geq \beta_n.$$
We say that $a$ {\it dominates} $b$ (denoted $a\trianglerighteq b$) 
if 
$$\sum_{i=1}^k \alpha_i \geq \sum_{i=1}^k \beta_i \quad \text{and} \quad k=1, \ldots, n-1 \quad 
\text{and} \quad \sum_{i=1}^n \alpha_i =\sum_{i=1}^n \beta_i.$$ Equivalently, $a \trianglerighteq b$ 
if $b$ is a convex combination of vectors obtained from $a$ by permutations of coordinates.

One can show that 
$$T\left(R_1, C_1; {\bold 1} \right) \geq T\left(R_2, C_2; {\bold 1} \right) 
\quad \text{provided} \quad R_2  \trianglerighteq R_1 \quad \text{and} \quad 
C_2 \trianglerighteq C_1. \tag2.4$$
The proof consists of two steps. First, assuming that 
$R=\left(r_1 \geq r_2 \geq \ldots \geq r_m\right)$ and $C=\left(c_1 \geq c_2 \geq \ldots \geq c_n \right)$
 we express $T(R, C; {\bold 1})$ in terms of 
{\it Kostka numbers},
$$T(R, C; {\bold 1} )= \sum_{A} K_{AR} K_{AC},$$
where the sum is taken over all $A=\left(a_1 \geq a_2 \geq \ldots \geq a_s \right)$, see
Section 6.I of \cite{Ma95}. Then we apply the inequality
$$K_{AB_2} \leq K_{AB_1} \quad \text{provided} \quad B_2  \trianglerighteq B_1,$$
see Section 7.I of \cite{Ma95}. Inequality (2.4) is consistent with the hypothesis (2.3).
\bigskip
To prove Theorem 2.1, we represent $T(R, C; W)$ as the expectation of the permanent 
of a random $N \times N$ matrix $A$ with exponentially distributed entries using a result from 
\cite{Ba05}. Then using the theory of matrix scaling \cite{MO68}, \cite{RS89}, \cite{L+00}, we represent 
$\per A$ as the product of a ``large and tame'' and a ``small and wild'' factors.
The ``tame'' factor contributes the bulk to the expectation and it satisfies the conditions of 
the Pr\'ekopa-Leindler inequality (1.1.2), the fact that ultimately results in the 
inequality of Theorem 2.1. The ``wild'' factor is harder to analyze, but  it  does not vary much since it lies within the 
low bound provided by the van der Waerden estimate \cite{Eg81}, \cite{Fa81} and the 
upper bound provided by the Bregman-Minc estimate \cite{Br73}. It contributes to the 
correction term in Theorem 2.1 and Corollary 2.2. 

We discuss preliminaries in Section 3 and present the proofs of Theorem 2.1 and Corollary 2.2
in Section 4.    

\head 3. A permanental representation of $T(R, C; W)$ \endhead

Recall that the {\it permanent} of an $N \times N$ matrix $A=\left(a_{ij} \right)$ is 
defined by 
$$\per A=\sum_{\sigma \in S_N} \prod_{i=1}^N a_{i \sigma(i)},$$
where $S_N$ is the symmetric group of all permutations of $\{1, \ldots, N\}$.
We say that a random variable $\gamma$ has the {\it standard exponential distribution} 
if 
$$\PP(\gamma > t)=\cases e^{-t} &\text{if\ } t>0 \\ 1 &\text{otherwise.} \endcases$$
The following result expressing $T(R, C; W)$ as the expectation of the permanent of 
a random matrix was proved in \cite{Ba05}.
\proclaim{(3.1) Theorem} Given a positive integer $m$-vector $R=(r_1, \ldots, r_m)$ and 
a positive integer $n$-vector $C=(c_1, \ldots, c_n)$ such that 
$$\sum_{i=1}^m r_i =\sum_{j=1}^n c_j =N,$$
and an $m \times n$ matrix $W=\left(w_{ij}\right)$,
we construct an $N \times N$ random matrix $A$ as follows: the set of rows of $A$ is represented 
as a disjoint 
union of $m$ subsets of cardinalities $r_1, \ldots, r_m$ whereas the set of columns of $A$ is 
represented as a disjoint union of $n$ subsets of cardinalities $c_1, \ldots, c_n$, so that 
$A$ is represented as a block matrix of $mn$ blocks $r_i \times c_j$. Let  
let $G=\left(g_{ij}\right)$ be the $m \times n$ matrix with $g_{ij}=w_{ij}\gamma_{ij}$, where 
$\gamma_{ij}$ are independent standard exponential random variables. We fill 
the $(i,j)$th block $r_i \times c_j$ of $A=A(G)$ by the copies of $g_{ij}$. 
Then 
$$T(R, C; W)={\EE \per A \over r_1! \cdots r_m! c_1! \cdots c_n!}.$$
\endproclaim

Next, we need some results on matrix scaling, in particular as described in \cite{MO68} and 
\cite{RS89}.
\subhead (3.2) Matrix scaling \endsubhead
Let $G=\left(g_{ij}\right)$ be a positive $m \times n$ matrix and 
let $r_1, \ldots, r_m$ and $c_1, \ldots, c_n$ be positive numbers such that 
$$\sum_{i=1}^m r_i =\sum_{j=1}^n c_j=N.$$
Then there exist a unique positive $m \times n$ matrix $L=\left(l_{ij}\right)$ and 
positive numbers $\mu_1, \ldots, \mu_m$ and $\lambda_1, \ldots, \lambda_n$ such that 
$$\split &\sum_{j=1}^n l_{ij}=r_i \quad \text{for} \quad i=1, \ldots, m, \\
&\sum_{i=1}^m l_{ij}= c_j \quad \text{for} \quad j=1, \ldots, n \endsplit$$
and such that 
$$g_{ij}=l_{ij} \mu_i \lambda_j \quad \text{for all} \quad i, j.$$

Moreover, the numbers $\lambda_i, \mu_j$ are unique up to a scaling 
$$\mu_i \longmapsto \mu_i \tau, \lambda_j \longmapsto \lambda_j \tau^{-1} 
\quad \text{for some} \quad \tau >0 \quad \text{and all} \quad i,j$$
and can be obtained as follows.

Let 
$$\split F(G; x, y)=\sum_{i=1}^m \sum_{j=1}^n &g_{ij} \xi_i \eta_j \quad \text{for} \\ 
&x=\left(\xi_1, \ldots, \xi_m\right)  \quad \text{and} \quad y=\left(\eta_1, \ldots, \eta_n\right).
\endsplit$$
Then $F(G; x, y)$ attains a unique minimum on the set of pairs $(x,y)$ of vectors defined by the equations
$$\prod_{i=1}^m \xi_i^{r_i}=1 \quad \text{and} \quad \prod_{j=1}^n \eta_j^{c_j}=1$$
and inequalities
$$\xi_i >0 \quad \text{for} \quad i=1, \ldots, m \quad 
\text{and} \quad \eta_j >0 \quad \text{for} \quad j=1, \ldots, n.$$
Assuming that $x^{\ast}=\left(\xi_1^{\ast}, \ldots, \xi_m^{\ast} \right)$ and 
$y^{\ast}=\left(\eta_1^{\ast}, \ldots, \eta_n^{\ast}\right)$ is the minimum point, we
may let 
$$\mu_i ={F\left(G; x^{\ast}, y^{\ast}\right) \over  N \xi_i} \quad \text{and} 
\quad \lambda_j ={1 \over \eta_j} \quad \text{for all} \quad i,j,$$
see \cite{RS89} and \cite{MO68}.

Finally, we need some estimates for permanents.
\subhead (3.3) Estimates for permanents \endsubhead
Recall that an $N \times N$ matrix $B=\left(b_{ij}\right)$ is called {\it doubly stochastic} if it is non-negative
$$b_{ij} \geq 0 \quad \text{for} \quad i,j=1, \ldots, N$$
and all row and column sums are equal to 1:
$$\split &\sum_{j=1}^N b_{ij}=1 \quad \text{for} \quad i=1, \ldots, N \quad \text{and} \\
&\sum_{i=1}^N b_{ij}=1 \quad \text{for} \quad j=1, \ldots, N. \endsplit$$
The van der Waerden conjecture proved by G.P. Egorychev \cite{Eg81} and 
D.I. Falikman \cite{Fa81} asserts that 
$$\per B \geq {N! \over N^N} \tag3.3.1$$
if $B$ is a doubly stochastic $N \times N$ matrix, see also Chapter 12 of \cite{LW01}.

The following upper bound was conjectured by H. Minc and proved by L.M. Bregman \cite{Br73}, 
see also Chapter 11 of \cite{LW01}. 

Let $B=\left(b_{ij}\right)$ be an $N \times N$ matrix such that $b_{ij} \in \{0, 1\}$ for all 
$i,j$ and let 
$$\sum_{j=1}^N b_{ij} =s_i \quad \text{for} \quad i=1, \ldots, N.$$
Then 
$$\per B \leq \prod_{i=1}^N \left(s_i! \right)^{1/s_i}.$$
We will need the following corollary of the Bregman-Minc inequality, see \cite{So03}.

Let $B=\left(b_{ij}\right)$ be an $N \times N$ matrix such that 
$$\split &\sum_{j=1}^N b_{ij}=1 \quad \text{for} \quad i=1, \ldots, N \quad \text{and} \\
&0 \leq b_{ij} \leq {1 \over s_i} \quad \text{for} \quad j=1, \ldots, N \endsplit$$
and positive integers $s_1, \ldots, s_N$.
Then 
$$\per B \leq \prod_{i=1}^N {(s_i!)^{1/s_i} \over s_i}. \tag3.3.2$$ 
Of course, similar estimates hold if we interchange rows and columns.

\head 4. Proofs \endhead

In this section, we prove Theorem 2.1 and Corollary 2.2.

\subhead (4.1) Notation \endsubhead
Given an $m \times n$ positive matrix $G=\left(g_{ij}\right)$, let us define
$$\split F(G; x, y)=\sum_{i=1}^m \sum_{j=1}^n g_{ij} &\xi_i \eta_j \quad 
\text{for} 
\\ x=&\left(\xi_1, \ldots, \xi_m\right) \quad \text{and} \quad 
y=\left(\eta_1, \ldots \eta_n \right). \endsplit$$
For positive vectors $R=\left(r_1, \ldots, r_m\right)$ and 
$C=\left(c_1, \ldots, c_n \right)$, we define 
$$ \split f(G; R, C)=&\min  F(G; x, y) \\
&\text{for} \quad x=\left(\xi_1, \ldots, \xi_m\right) \quad \text{and} \quad y=\left(\eta_1, \ldots, \eta_n \right) \\
&\text{subject to} \quad \prod_{i=1}^m \xi_i^{r_i}=\prod_{j=1}^n \eta_j^{c_j} =1, \endsplit$$
see Section 3.2.
We recall notation
$$|R|=\sum_{i=1}^m r_i \quad \text{and} \quad |C|=\sum_{j=1}^n c_j.$$

First, we establish a certain convexity property of $f(G; R, C)$.

\proclaim{(4.2) Lemma} Let $G_1, \ldots, G_p$ be positive $m \times n$ matrices,
let $R_1, \ldots, R_p$ be positive $m$-vectors, and let $C_1, \ldots, C_p$ be positive 
$n$-vectors such that 
$$|R_1| =\ldots = |R_p| =|C_1| = \ldots = |C_p|.$$
Suppose further that $\alpha_1, \ldots, \alpha_p \geq 0$ are numbers such that 
$\alpha_1 + \ldots + \alpha_p=1$.
Let us define
$$G=\sum_{k=1}^p \alpha_k G_k, \quad R=\sum_{k=1}^p \alpha_k R_k, \quad \text{and} \quad 
C=\sum_{k=1}^p \alpha_k C_k.$$
Then 
$$f(G; R, C) \geq \prod_{k=1}^p f^{\alpha_k}\left(G_k; R_k, C_k\right).$$
\endproclaim
\demo{Proof}
Suppose that 
$$\split &R_k=\left(r_{1k}, \ldots, r_{mk} \right), \quad 
C_k=\left( c_{1k}, \ldots, c_{nk} \right), \quad 
R=\left(r_1, \ldots, r_m\right), \quad \text{and} \\ 
&C=\left(c_1, \ldots, c_n \right). \endsplit$$ 
In particular,
$$\aligned &r_i=\sum_{k=1}^p \alpha_k r_{ik} \quad \text{for} \quad i=1, \ldots, m 
\quad \text{and} \\ &c_j=\sum_{k=1}^p \alpha_k c_{jk} \quad \text{for} \quad 
j=1, \ldots, n. \endaligned \tag4.2.1 $$

 Let $x=\left(\xi_1, \ldots, \xi_m\right)$ and $y=\left(\eta_1, \ldots, \eta_n\right)$ be positive vectors 
 such that 
 $$\prod_{i=1}^m \xi_i^{r_i}= \prod_{j=1}^n \eta_j^{c_j}=1. \tag4.2.2$$
Then 
$$F(G; x, y)=\sum_{k=1}^p \alpha_k F(G_k; x, y) \geq \prod_{k=1}^p 
F^{\alpha_k}(G_k; x, y).$$  
Let 
$$t_k=\left( \prod_{i=1}^m \xi_i^{r_{ik}} \right)^{1/|R|} \quad \text{and} \quad 
s_k= \left( \prod_{j=1}^n \eta_j^{c_{jk}} \right)^{1/|C|} \quad \text{for} \quad k=1, \ldots, p.$$ 
Then 
$$F(G_k; x, y) = t_k s_k F\left(G_k; t_k^{-1} x, s_k^{-1}y\right) \geq 
t_k s_k f\left(G_k; R_k, C_k \right),$$
since vectors $t_k^{-1} x$ and $s_k^{-1}y$ satisfy (4.2.2) with $r_i$ and $c_j$ replaced 
by $r_{ik}$ and $c_{jk}$ respectively.
Therefore,
$$F(G; x,y) \geq \prod_{k=1}^p t_k^{\alpha_k} s_k^{\alpha_k} f^{\alpha_k}\left(G_k; R_k, C_k \right).$$ 
On the other hand, by (4.2.1) and (4.2.2), we have 
$$\prod_{k=1}^p t_k^{\alpha_k} =\left(\prod_{i=1}^m \xi_i^{\sum_{k=1}^p \alpha_k r_{ik}}\right)^{1/|R|}
=\left(\prod_{i=1}^m \xi_i^{r_i}\right)^{1/|R|}=1,$$
and, similarly, 
 $$\prod_{k=1}^p s_k^{\alpha_k} =\left(\prod_{j=1}^n \eta_j^{\sum_{k=1}^p \alpha_k c_{jk}}\right)^{1/|C|}
=\left(\prod_{j=1}^n \eta_j^{c_j}\right)^{1/|C|}=1.$$

Since the inequality 
$$F(G; x,y) \geq \prod_{k=1}^p f^{\alpha_k}\left(G_k, R_k, C_k \right)$$
holds for any positive $x$ and $y$ satisfying (4.2.2), the proof follows.
{\hfill \hfill \hfill} \qed
\enddemo

Next, we consider block matrices $A$ as in Theorem 3.1.

\proclaim{(4.3) Lemma} Let $G=\left(g_{ij} \right)$ be an $m \times n$ positive matrix.
Let $R=\left(r_1, \ldots, r_m \right)$ and $C=\left(c_1, \ldots, c_n \right)$ be positive 
integer vectors such that $|R|=|C|=N$. Let us consider the $N \times N$ block matrix $A$, where 
the $(i,j)$th block of size $r_i \times c_j$ is filled by copies of $g_{ij}$. Then there exists an 
$N \times N$ block matrix $B$ with the same block structure as $A$ and such that
\roster
\item Matrix $B$ is doubly stochastic;
\item The entries in the $(i,j)$th block of $B$ do not exceed $\min\{1/r_i,\ 1/c_j\}$;
\item We have 
$$\per A= N^{-N} f^N(G; R, C) \left( \prod_{i=1}^m r_i^{r_i} \right) \left(\prod_{j=1}^n c_j^{c_j}\right)
\per B$$
\item 
$${N! \over N^N} \leq \per B \leq \min \left\{\prod_{i=1}^m {r_i! \over r_i^{r_i}},  \quad \prod_{j=1}^n {c_j! \over c_j^{c_j}} \right\}.$$ 
\endroster
\endproclaim

\demo{Proof}  Let $L=\left(l_{ij}\right)$ be the $m \times n$ matrix and let 
$\mu_i$, $i=1, \ldots, m$, and $\lambda_j$, $j=1, \ldots, n$, be numbers such that
$$g_{ij}=l_{ij} \mu_i \lambda_j \quad \text{for} \quad i=1, \ldots, m \quad \text{and} \quad 
j=1, \ldots, n$$ 
and such that
$$\split &\sum_{j=1}^n l_{ij}=r_i \quad \text{for} \quad i=1, \ldots, m \quad \text{and} \\
&\sum_{i=1}^m l_{ij}=c_j \quad \text{for} \quad j=1, \ldots, n, \endsplit$$
see Section 3.2.

Let us divide the entries in the $(i,j)$th block of $A$ by the product $\mu_i r_i \lambda_j c_j$. We get 
the matrix $B$ with the entries in the $(i,j)$th block equal to $l_{ij}/r_i c_j$. It is seen now that $B$ is 
doubly stochastic and that the entries in the $(i,j)$th block of $B$ do not exceed 
$\min\{1/r_i,\ 1/c_j\}$. Furthermore,
$$\per A= \left(\prod_{i=1}^m \left(\mu_i r_i\right)^{r_i} \right) 
\left(\prod_{j=1}^n \left(\lambda_j c_j\right)^{c_j} \right) \per B.$$ 
On the other hand, if one computes $\mu_i$ and $\lambda_j$ by optimizing 
$F(G; x, y)$  as in Section 3.2, one gets 
$$\prod_{i=1}^m \mu_i^{r_i} ={f^N(G; R, C) \over N^N}  \quad \text{and} \quad 
\prod_{j=1}^n \lambda_j^{c_j}=1,$$
which completes the proof of Part (3).

Part (4) follows by Parts (1) and (2) and estimates (3.3.1) and (3.3.2).
{\hfill \hfill \hfill} \qed
\enddemo

Now we are ready to prove Theorem 2.1.

\demo{Proof of Theorem 2.1} Without loss of generality, we assume that 
$w_{ij}>0$ for all $i,j$. 

In the space $\Mat(m,n)$ of $m \times n$ real matrices $G=\left(g_{ij}\right)$, we consider the exponential measure $dG$ with the 
density 
$$\prod_{ij} w_{ij}^{-1} \exp\left\{-g_{ij}/w_{ij} \right\} \quad \text{if} \quad g_{ij} >0 \quad \text{for all} \quad i,j$$
and 0 elsewhere.

We note that $dG$ is a log-concave measure.

Given positive integer vectors $R=\left(r_1, \ldots, r_m \right)$ and $C=\left(c_1, \ldots, c_n \right)$
and a positive $m \times n$ matrix $G$, let  $A(G; R, C)$ be the $N \times N$ block matrix constructed 
as in Theorem 3.1. Then, by Theorem 3.1,
$$T(R, C; W)=\left(\prod_{i=1}^m {1 \over r_i!} \right)
\left( \prod_{j=1}^n {1 \over c_j!} \right) \int_{\Mat(m,n)} \per A(G; R, C) \ d G.$$
From Lemma 4.3,
$$\split T(R, C; W) \geq &{N! \over N^N} N^{-N} \left(\prod_{i=1}^m {r_i^{r_i} \over r_i!} \right) 
\left(\prod_{j=1}^n {c_j^{c_j} \over c_j!} \right) \\
&\times  \int_{\Mat(m,n)} f^N(G; R, C) \ d G \\
= &{N! \over N^N} N^{-N}  \omega(R) \omega(C) \int_{\Mat(m,n)} f^N(G; R, C) \ d G. \endsplit $$

Similarly,
letting $R_k=\left(r_{1k}, \ldots, r_{mk} \right)$ and 
$C_k=\left( c_{1k}, \ldots, c_{nk} \right)$, by Theorem 3.1 we obtain 
$$T(R_k, C_k; W_k)=\left(\prod_{i=1}^m {1 \over r_{ik}!} \right)
\left( \prod_{j=1}^n {1 \over c_{jk}!} \right) \int_{\Mat(m,n)} \per A(G; R_k, C_k) \ d G,$$
and from Lemma 4.3
$$\split T\left(R_k, C_k; W \right) \leq &N^{-N} \left(\prod_{i=1}^m {r_{ik}^{r_{ik}} \over r_{ik}!} \right)
\left( \prod_{j=1}^n {c_{jk}^{c_{jk}} \over c_{jk}!} \right) \\
&\times \min \left\{ \prod_{i=1}^m {r_{ik}! \over r_{ik}^{r_{ik}}}, \quad \prod_{j=1}^n {c_{jk}! \over c_{jk}^{c_{jk}}} \right\} 
\int_{\Mat(m,n)} f^N(G; R_k, C_k) \ d G\\
= & N^{-N} \min\bigl\{ \omega(R_k),\  \omega(C_k) \bigr\}  
\int_{\Mat(m,n)} f^N(G; R_k, C_k) \ d G. \endsplit$$

By Lemma 4.2, for any positive matrices $G_1, \ldots, G_k$ we have 
$$f(G; R, C) \geq \prod_{k=1}^p f^{\alpha_k} \left(G_k; R_k, C_k \right), \quad 
\text{where} \quad G=\sum_{k=1}^p \alpha_k G_k.$$

Applying the Pr\'ekopa-Leindler inequality (1.1.2), we obtain
$$\int_{\Mat(m,n)} f^N(G; R, C) \ d G \geq 
\prod_{k=1}^p \left( \int_{\Mat(m,n)} f^N \left(G; R_k, C_k\right) d G \right)^{\alpha_k}.$$ 
Therefore,
$${N^N  \over N!} {T(R, C; W) \over \omega(R) \omega(C)} \geq 
\prod_{k=1}^p \left( {T\left(R_k, C_k; W \right) \over \min\bigl\{ \omega(R_k),\  \omega(C_k) \bigr\} }\right)^{\alpha_k}$$ 
and the proof follows.
{\hfill \hfill \hfill} \qed
\enddemo

\demo{Proof of Corollary 2.2} We use that the function 
$$s \longmapsto {b^b \over \Gamma(b+1)}, \quad b >0$$ 
is log-convex. Therefore, the function
$$\omega\left(b_1, \ldots, b_l\right)= 
\prod_{i=1}^l {b_i^{b_i} \over \Gamma(b_i+1)}$$
is log-convex on the positive orthant $b_1>0, \ldots, b_l>0$.

Thus we have
$$\omega(R) \geq \prod_{k=1}^p \omega^{\alpha_k}\left(R_k\right) \quad \text{and} \quad 
\omega(C) \geq \prod_{k=1}^p \omega^{\alpha_k}\left(C_k\right).$$
Hence
$$
\split &{N^N \over N!} {T(R,C; W) \over \omega(C)} \geq \prod_{k=1}^p T^{\alpha_k}\left(R_k, C_k; W \right) \quad \text{and} \\
&{N^N \over N!} {T(R,C; W) \over \omega(R)} \geq \prod_{k=1}^p T^{\alpha_k}\left(R_k, C_k; W \right),
\endsplit$$
from which Part (1) follows.

Similarly, since $\omega$ is log-convex,
$$\omega(R)  \geq \omega\left( |R|/m, \ldots, |R|/m \right) \quad \text{and} \quad 
\omega(C) \geq \omega\left(|C|/n, \ldots, |C|/n \right),$$
from which Part (2) follows.

Finally, by Stirling's formula
$$(2 \pi s)^{1/2} s^s e^{-s} e^{1 \over 12s+1} < \Gamma(s+1) < (2 \pi s)^{1/2} s^s e^{-s} e^{1 \over 12s}$$
and Part(3) follows.
{\hfill \hfill \hfill} \qed
\enddemo

\head Acknowledgment \endhead

I am grateful to Alexander Yong for many useful conversations.


\Refs
\widestnumber\key{AAAA}

\ref\key{Ba05}
\by A. Barvinok
\paper Enumerating contingency tables via random permanents
\paperinfo preprint arXiv math.CO/0511596
\yr 2005
\endref

\ref\key{Be74}
\by E.A. Bender
\paper The asymptotic number of non-negative integer matrices with given row and column sums 
\jour Discrete Math. 
\vol 10 
\yr 1974
\pages  217--223
\endref

\ref\key{Br73}
\by L.M. Bregman 
\paper Certain properties of
nonnegative matrices and their permanents 
\jour Dokl. Akad. Nauk
SSSR
\vol 211
 \yr 1973
 \pages 27--30
 \endref

\ref\key{B+04}
\by W. Baldoni-Silva, J.A. De Loera, and M. Vergne 
\paper Counting integer flows in networks 
\jour Found. Comput. Math. 
\vol 4 
\yr 2004
\pages 277--314
\endref

 \ref\key{DE85}
 \by P. Diaconis and B. Efron
 \paper Testing for independence in a two-way table: new interpretations of the chi-square statistic. With discussions and with a reply by the authors
 \jour  Ann. Statist. 
 \vol 13 
 \yr 1985
 \pages  845--913
\endref

\ref\key{DG95}
\by P. Diaconis and A. Gangolli
\paper Rectangular arrays with fixed margins
\inbook Discrete probability and algorithms (Minneapolis, MN, 1993)
\pages 15--41 
\bookinfo IMA Vol. Math. Appl.
\vol  72 
\publ Springer
\publaddr  New York
\yr  1995
\endref 

\ref\key{DG04}
\by P. Diaconis and A. Gamburd 
\paper Random matrices, magic squares and matching polynomials 
\jour Electron. J. Combin. 
\vol 11 
\yr 2004
\paperinfo Research Paper 2
\pages 26 pp
\endref

\ref\key{Eg81} 
\by G.P. Egorychev 
\paper The solution of van der
Waerden's problem for permanents 
\jour Adv. in Math. 
\vol 42 
\yr 1981 
\pages 299--305 
\endref

\ref\key{Fa81}
 \by D.I. Falikman 
 \paper Proof of the van der
Waerden conjecture on the permanent of a doubly stochastic matrix
(Russian)
 \jour Mat. Zametki
 \vol 29 
 \yr 1981
  \pages 931--938
\endref

\ref\key{Ga02}
\by R.J.  Gardner 
\paper The Brunn-Minkowski inequality
\jour Bull. Amer. Math. Soc. (N.S.) 
\vol 39 
\yr 2002
\pages  355--405
\endref 

\ref\key{GG01}
\by R.J. Gardner and P. Gronchi 
\paper A Brunn-Minkowski inequality for the integer lattice 
\jour Trans. Amer. Math. Soc. 
\vol 353 
\yr 2001
\pages 3995--4024
\endref

\ref\key{Le01}
\by M. Ledoux 
\book The Concentration of Measure Phenomenon 
\bookinfo Mathematical Surveys and Monographs
\vol  89 
\publ American Mathematical Society
\publaddr  Providence, RI
\yr  2001
\endref

\ref\key{LW01}
 \by J.H. van Lint and R.M. Wilson 
 \book A Course in
Combinatorics. Second edition 
\publ Cambridge University Press
\publaddr Cambridge
 \yr 2001 
 \endref

\ref\key{L+00} 
\by N. Linial, A. Samorodnitsky, and A. Wigderson
\paper A deterministic strongly polynomial algorithm for matrix
scaling and approximate permanents
 \jour Combinatorica
 \vol 20
 \yr 2000 
 \pages 545--568 
 \endref

\ref\key{Ma95}
\by I.G. Macdonald
\book Symmetric Functions and Hall Polynomials. Second edition. With contributions by A. Zelevinsky
\bookinfo Oxford Mathematical Monographs. Oxford Science Publications
\publ The Clarendon Press, Oxford University Press
\publaddr New York
\yr 1995
\endref

\ref\key{MO68}
\by A. Marshall and I. Olkin
\paper Scaling of matrices to achieve specified row and column sums 
\jour Numer. Math. 
\vol 12 
\yr 1968 
\pages 83--90
\endref

\ref\key{RS89}
\by U. Rothblum and H. Schneider
\paper Scalings of matrices which have prespecified row sums and column sums via optimization
\jour Linear Algebra Appl. 
\vol 114/115 
\yr 1989
\pages 737--764
\endref

\ref\key{St97} 
\by R.P. Stanley
\book Enumerative Combinatorics. Vol. 1
\bookinfo Cambridge Studies in Advanced Mathematics
\vol 49
\publ Cambridge University Press
\publaddr Cambridge
\yr 1997
\endref

\ref\key{So03} 
\by G.W. Soules 
\paper New permanental upper bounds
for nonnegative matrices
\jour Linear Multilinear Algebra
\vol 51
\yr 2003 
\pages 319--337 
\endref


\ref\key{Vi03}
\by C. Villani 
\book Topics in Optimal Transportation 
\bookinfo Graduate Studies in Mathematics
\vol 58 
\publ American Mathematical Society
\publaddr Providence, RI
\yr  2003
\endref

\endRefs
\enddocument
\end